\documentclass[12pt,reqno]{amsart}

\usepackage{amssymb}

\usepackage{mathrsfs}

\usepackage{bbm}

\DeclareMathAlphabet{\mathpzc}{OT1}{pzc}{m}{it}

\def\fm{{\mathfrak m}}

\def\fp{{\mathfrak p}}

\def\amp{\operatorname{amp}}

\def\b{\operatorname{b}}

\def\cmd{\operatorname{cmd}}

\def\D{\operatorname{D}}
\def\Dc{\D^{\operatorname{c}}}

\def\depth{\operatorname{depth}}
\def\Dfg{\D^{\operatorname{fg}}}
\def\Dfgm{\Dfg_{-}}
\def\Dfgp{\Dfg_{+}}
\def\Dfgpm{\Dfg_{\b}}
\def\dim{\operatorname{Kdim}}

\def\H{\operatorname{H}}

\def\Hom{\operatorname{Hom}}

\def\inf{\operatorname{inf}}

\def\kid{\operatorname{k.\!id}}
\def\kpd{\operatorname{k.\!pd}}

\def\LTensor{\stackrel{\operatorname{L}}{\otimes}}
\def\opp{\operatorname{o}}
\def\pd{\operatorname{pd}}

\def\RHom{\operatorname{RHom}}

\def\sup{\operatorname{sup}}
\def\Supp{\operatorname{Supp}}

\numberwithin{equation}{part}

\newtheorem{Lemma}{Lemma}[section]
\newtheorem{Theorem}[Lemma]{Theorem}
\newtheorem{Proposition}[Lemma]{Proposition}

\theoremstyle{definition}

\newtheorem{Setup}[Lemma]{Setup}

\newtheorem{Remark}[Lemma]{Remark}

\newtheorem{Notation}[Lemma]{Notation}

\def\rk{A}

\begin{document}

\setlength{\parskip}{6pt}

\title[Amplitude inequalities]
{Amplitude inequalities for Differential Graded modules}

\author{Peter J\o rgensen}
\address{Department of Pure Mathematics, University of Leeds,
Leeds LS2 9JT, United Kingdom}
\email{popjoerg@maths.leeds.ac.uk}
\urladdr{http://www.maths.leeds.ac.uk/\~{ }popjoerg}




\begin{abstract} 

Differential Graded Algebras can be studied through their Differential
Graded modules.  Among these, the compact ones attract particular
attention.  

This paper proves that over a suitable chain Differential Graded
Algebra $R$, each compact Differential Graded module $M$ satisfies
$\amp M \geq \amp R$, where $\amp$ denotes amplitude which is
defined in a straightforward way in terms of the homology of a
DG module.

In other words, the homology of each compact DG module $M$ is at least
as long as the homology of $R$ itself.  Conversely, DG modules with
shorter homology than $R$ are not compact, and so in general, there
exist DG modules with finitely generated homology which are not
compact.

Hence, in contrast to ring theory, it makes no sense to define finite
global dimension of DGAs by the condition that each DG module with
finitely generated homology must be compact.

\end{abstract}

\maketitle

\setcounter{section}{-1}
\section{Introduction}
\label{sec:introduction}

\noindent 

Differential Graded Algebras (DGAs) play an important role in both
ring theory and algebraic topology.  For instance, if $M$ is a complex
of modules, then the endomorphism complex $\Hom(M,M)$ is a DGA with
multiplication given by composition of endomorphisms, and this can be
used to prove ring theoretical results, see \cite{KellerTilting} and
\cite{KellerDG}.  Another example is that over a commutative ring, the
Koszul complex on a series of elements is a DGA, see \cite[sec.\
4.5]{Weibel}, and again, ring theoretical results ensue, see
\cite{Iyengar}.

Likewise, DGAs occur naturally in algebraic topology, where the
canonical example is the singular cochain complex $C^*(X)$ of a
topological space $X$.  Other constructions also give DGAs; for
instance, if $G$ is a topological monoid, then the singular chain
complex $C_*(G)$ is a DGA whose multiplication is induced by the
composition of $G$; see \cite{FHTbook}.

Just as rings can be studied through their modules, DGAs can be
studied through their Differential Graded modules (DG modules), and
this is the subject of the present paper.  

The main results are a number of ``amplitude inequalities'' which give
bounds on the amplitudes of various types of DG modules.  Such results
have been known for complexes of modules over rings since Iversen's
paper \cite{Iversen}, and it is natural to seek to extend them to DG
modules.

Another main point, implied by one of the amplitude inequalities, is
that, in contrast to ring theory, it appears to make no sense to
define finite global dimension of DGAs by the condition that each DG
module with finitely generated homology must be compact.  This is of
interest since several people have been asking how one might define
finite global dimension for DGAs.

\medskip
\noindent
{\em First main Theorem. }
To get to the first main Theorem of the paper, recall from
\cite{KellerDG} that if $R$ is a DGA then a good setting for DG
modules over $R$ is the derived category of DG left-$R$-modules
$\D(R)$.

A DG left-$R$-module is called compact if it is in the smallest
triangulated subcategory of $\D(R)$ containing $R$, or, to use the
language of topologists, if it can be finitely built from $R$.  The
compact DG left-$R$-modules form a triangulated subcategory $\Dc(R)$
of $\D(R)$, and play the same important role as finitely presented
modules of finite projective dimension do in ring theory.

The amplitude of a DG module $M$ is defined in terms
of the homology $\H(M)$ by
\[
  \amp M = \sup \{\, i \mid \H_i(M) \not= 0 \,\}
           - \inf \{\, i \mid \H_i(M) \not= 0 \,\}.
\]

\noindent
{\bf Theorem A. }
{\em 
Assume $\amp R < \infty$.  Let $L$ be in $\Dc(R)$ and suppose $L \not\cong
0$.  Then 
\[
  \amp L \geq \amp R.
\]
}

\noindent
Expressed in words, this says that among the compact DG modules, none
can be shorter than $R$ itself.  The Theorem will be proved in the
situation specified in Setup \ref{set:blanket} below; the main point
is that $R$ is a local DGA, that is, a chain DGA for which $\H_0(R)$
is a local commutative noetherian ring.  The multiplication in
$\H_0(R)$ is induced by the multiplication in $R$.

Of equal significance to Theorem A is perhaps the following
consequence: If $\amp R \geq 1$, that is, if $R$ is a true DGA in the
sense that it is not quasi-isomorphic to a ring, then a DG module with
amplitude zero cannot be compact.  There are many such DG modules and
they can even be chosen so that their homology $\H(M)$ is
finitely generated as a module over the ring $\H_0(R)$.  The scalar
multiplication of $\H_0(R)$ on $\H(M)$ is induced by the scalar
multiplication of $R$ on $M$.  For a concrete example, note that
$\H_0(R)$ itself can be viewed as a DG module via the canonical
surjection $R \rightarrow \H_0(R)$ which exists because $R$ is a chain DGA.

So if $\amp R \geq 1$ then there are DG modules with finitely
generated homology over $\H_0(R)$ which are not compact.  Hence, as
mentioned above, it appears to make no sense to define finite global
dimension of DGAs by the condition that each DG module with finitely
generated homology must be compact.  This contrasts sharply with ring
theory where this precise definition works for several classes of
rings, such as the local commutative noetherian ones.

\medskip
\noindent
{\em Second main Theorem. }
To explain the second main Theorem of the paper, let me first give an
alternative, equivalent formulation of Theorem A.  

Let $\Dfgp(R)$ denote the DG left-$R$-modules $M$ for which each
$\H_i(M)$ is finitely generated over $\H_0(R)$, and for which $\H_i(M)
= 0$ for $i \ll 0$.  It turns out that the compact DG left-$R$-modules
are exactly the DG modules in $\Dfgp(R)$ for which $\kpd_R M <
\infty$, where $\kpd$ denotes $k$-projective dimension, see Notation
\ref{not:blanket} and Lemma \ref{lem:kpd}.  Hence the following is an
equivalent formulation of Theorem A.

\noindent
{\bf Theorem A${}^{\prime}$. }
{\em 
Assume $\amp R < \infty$.  Let $L$ be in $\Dfgp(R)$ and suppose $\kpd_R
L < \infty$ and $L \not\cong 0$.  Then
\[
  \amp L \geq \amp R.
\]
}

\noindent
The dual of Theorem A${}^{\prime}$ is now the second main Theorem of
the paper, which will also be proved in the situation specified in
Setup \ref{set:blanket}.  To state it, some more notation is
necessary. 

Let $\rk$ denote a local commutative noetherian ground ring over which
$R$ is a DGA.  Let $\Dfgm(R)$ be the DG left-$R$-modules $M$ for which
each $\H_i(M)$ is finitely generated over $\H_0(R)$, and for which
$\H_i(M) = 0$ for $i \gg 0$.  Let $\kid$ denote $k$-injective
dimension, see Notation \ref{not:blanket}, and let $\cmd_{\rk}\!R$
denote the so-called Cohen-Macaulay defect of $R$ over $\rk$, see
\cite[(1.1)]{AvrFox} or Notation \ref{not:dim}.

\noindent
{\bf Theorem B. }
{\em 
Assume $\amp R < \infty$.  Let $I$ be in $\Dfgm(R)$ and suppose
$\kid_R I < \infty$ and $I \not\cong 0$.  Then
\[
  \amp I \geq \cmd_{\rk}R.
\]
}

\medskip
\noindent
{\em Comments and connections. }
Theorems A, A${}^{\prime}$, and B will be obtained as corollaries of a
more general amplitude inequality, Theorem \ref{thm:amp}, which is a
DGA generalization of the first of Iversen's amplitude inequalities
from \cite{Iversen}.

Theorem B can be written in a more evocative form for certain DGAs.
Suppose that the ground ring $\rk$ has a dualizing complex $C$ and
consider $D = \RHom_{\rk}(R,C)$ which is sometimes a so-called
dualizing DG module for $R$, see \cite{FIJ}.  Since $\amp R < \infty$
implies $\cmd_{\rk}R = \amp \RHom_{\rk}(R,C) = \amp D$ by
\cite[(1.3.2)]{AvrFox}, Theorem B takes the form
\[
  \amp I \geq \amp D.
\]
So if $D$ is indeed a dualizing DG module for $R$, then, expressed in
words, Theorem B says that among the DG modules in $\Dfgm(R)$ with
$\kid_R I < \infty$, none can be shorter than the dualizing DG module
$D$.

Theorems A, A$^{\prime}$, and B are specific to the finite amplitude
case, and fail completely if the amplitude of $R$ is permitted to be
infinite: Let $K$ be a field and consider the polynomial algebra
$K[X]$ as a DGA where $X$ is placed in homological degree $1$ and
where the differential is zero.  There is a distinguished triangle
\[
  \Sigma K[X] \rightarrow K[X] \rightarrow K \rightarrow
\]
in $\D(K[X])$, involving $K[X]$, the suspension $\Sigma K[X]$ and the
trivial DG module $K$.  Applying the functor $\RHom_{K[X]}(-,K)$ gives
a distinguished triangle
\[
  \RHom_{K[X]}(K,K) \rightarrow K \rightarrow \Sigma^{-1}K \rightarrow
\]
which shows that $\RHom_{K[X]}(K,K)$ has bounded homology whence
$\kpd_{K[X]}K < \infty$ and $\kid_{K[X]}K < \infty$, see Notation
\ref{not:blanket}. 

However, it is clear that $\amp K = 0$ and $\amp K[X] = \infty$, and
not hard to show $\cmd_{K}K[X] = \infty$.  Hence
\[
  \amp K < \amp K[X],
\]
showing that Theorems A${}^{\prime}$ and A fail, and
\[
  \amp K < \cmd_K K[X],
\]
showing that Theorem B fails.

Note that while Theorems A and A${}^{\prime}$ are uninteresting if $R$
is a ring concentrated in degree zero (for which $\amp R = 0$),
Theorem B is already interesting in this case.  For instance, if $R$
is equal to the ground ring $\rk$ placed in degree zero, then Theorem
B says that $\amp I \geq \cmd \rk$ when $I \not\cong 0$ is a complex
in $\Dfgm(\rk)$ with $\kid_{\rk} I < \infty$.

This implies the classical conjecture by Bass that if $\rk$ has a
finitely generated module $M$ with finite injective dimension, then
$\rk$ is a Cohen-Macaulay ring.  To see so, apply Theorem B to the
injective resolution $I$ of $M$.  This gives $0 = \amp I \geq \cmd
\rk$ whence $\cmd \rk = 0$, that is, $\rk$ is Cohen-Macaulay.  The
ring case of Theorem B and the fact that it implies the conjecture by
Bass has been known for a good while to commutative ring theorists,
but a published source seems hard to find.

The paper is organized as follows.  Section \ref{sec:notation}
explains some of the notation and terminology.  Sections
\ref{sec:homological_estimates} and \ref{sec:flat_base_change} are
preparatory; they prove a number of homological estimates and set up
some base change machinery.  Finally, Section \ref{sec:amp_ineq_A}
proves Theorems A and A${}^{\prime}$, and Section \ref{sec:amp_ineq_B}
proves Theorem B.

\medskip
\noindent
{\em Acknowledgement. }  I thank Henrik Holm for answering a question
and correcting some misprints, and Bernhard Keller for a conversation
about this material.

\section{Background}
\label{sec:notation}

This Section explains some of the notation and terminology of the
paper.  The usage will be standard and corresponds largely to such
references as \cite{Aldrich}, \cite{FIJ}, \cite{FJiic}, and
\cite{KellerDG}.

All proofs will be given under the following Setup.  Note, however,
that the results also hold in some other situations, see Remark
\ref{rmk:generalization}.

\begin{Setup}
\label{set:blanket}
By $\rk$ is denoted a local commutative noetherian ring, and by $R$ a
chain DGA (that is, $R_i = 0$ for $i \leq -1$) over $\rk$ for which
each $\H_i(R)$ is a finitely generated $\rk$-module.

It will be assumed that the canonical ring homomorphism $\rk
\rightarrow \H_0(R)$ is surjective.
\hfill $\Box$
\end{Setup}

\begin{Remark}
\label{rmk:A-structure}
Since $\rk$ is the ground ring for $R$, everything in sight will have
an $A$-structure.
\hfill $\Box$
\end{Remark}

\begin{Remark}
Since $\rk \rightarrow \H_0(R)$ is surjective, $\H_0(R)$ is a quotient
of $\rk$ and so $\H_0(R)$ is a local commutative noetherian ring.

In fact, $\rk \rightarrow \H_0(R)$ is equal to the composition $\rk
\rightarrow R \rightarrow \H_0(R)$, where $R \rightarrow \H_0(R)$ is
the canonical surjection which exists because $R$ is a chain DGA.
Through these morphisms, the residue class field $k$ of the local ring
$\H_0(R)$ can be viewed as a DG bi-module over $R$, and as a module
over $\rk$, and these will be denoted simply by $k$.

Note that $k$ viewed as a $\rk$-module is the residue class field of
$\rk$, because $\rk \rightarrow \H_0(R)$ is surjective.
\hfill $\Box$
\end{Remark}

\begin{Remark}
Any ring can be viewed as a trivial DGA concentrated in
degree zero.  A DG module over such a DGA is just a complex of
modules over the ring in question.

Also, an ordinary module over a ring can be viewed as a complex
concentrated in degree zero, and hence as a DG module over the ring
viewed as a DGA.

When there is more than one ring or DGA action on the same object, I
will sometimes use subscripts to indicate the actions.  For instance,
a DG left-$R$-right-$R$-module might be denoted ${}_{R}D_{R}$.
\hfill $\Box$
\end{Remark}

\begin{Notation}
\label{not:blanket}
By $R^{\opp}$ is denoted the opposite DGA of $R$ with product defined
in terms of the product of $R$ by $r \cdot s = (-1)^{\mid r \mid \mid
s \mid}sr$.  DG left-$R^{\opp}$-modules will be identified with DG
right-$R$-modules, and $\D(R^{\opp})$, the derived category of DG
left-$R^{\opp}$-modules, will be identified with the derived category
of DG right-$R$-modules.


The supremum and the infimum of the DG module $M$ are defined by
\begin{eqnarray*}
  \sup M & = & \sup \{\, i \mid \H_i(M) \not= 0 \,\}, \\
  \inf M & = & \inf \{\, i \mid \H_i(M) \not= 0 \,\};
\end{eqnarray*}
in these terms the amplitude of $M$ is
\[
  \amp M = \sup M - \inf M.
\]

The $k$-projective dimension, the $k$-injective dimension, and the
depth are defined by
\begin{eqnarray*}
  \kpd_R M & = & - \inf \RHom_R(M,k), \\
  \kid_R M & = & - \inf \RHom_R(k,M), \\
  \depth_R M & = & - \sup \RHom_R(k,M),
\end{eqnarray*}
see \cite[def.\ (1.1)]{FJiic}.  Here $\RHom$ is the right derived
functor of $\Hom$ which will be used along with $\LTensor$, the left
derived functor of $\otimes$.

Derived functors are defined on derived categories.  Some notation for
these was already given above, but let me collect it here.  The
derived category of DG left-$R$-modules is denoted by $\D(R)$.  The
full subcategory of compact objects is denoted by $\Dc(R)$.

The full subcategory of $\D(R)$ consisting of DG modules with each
$\H_i(M)$ finitely generated over $\H_0(R)$ and $\H_i(M) = 0$ for $i
\ll 0$ is denoted by $\Dfgp(R)$.  The full subcategory of $\D(R)$
consisting of DG modules with each $\H_i(M)$ finitely generated over
$\H_0(R)$ and $\H_i(M) = 0$ for $i \gg 0$ is denoted by $\Dfgm(R)$.
Finally, $\Dfgm(R) \cap \Dfgp(R)$ is denoted by $\Dfgpm(R)$.

If DG modules are viewed as having the differentials pointing to the
right, then $\Dfgp(R)$ consists of DG modules with homology extending
to the left, $\Dfgm(R)$ consists of DG modules with homology extending
to the right, and $\Dfgpm(R)$ consists of DG modules with bounded
homology.

Observe that $R$ could just be $\rk$ or $\H_0(R)$ concentrated in
degree zero.  Hence the notations introduced so far can also be
applied to $\rk$ and $\H_0(R)$, and define triangulated subcategories
$\Dc$, $\Dfgp$, $\Dfgm$, $\Dfgpm$ of the derived categories $\D(\rk)$
and $\D(\H_0(R))$.
\hfill $\Box$
\end{Notation}

\begin{Notation}
\label{not:dim}
It is well known that homological invariants such as projective
dimension (often denoted $\pd$) and depth can be extended from modules
to complexes of modules, see for instance \cite{Foxby} and
\cite{Iversen}.

I will need two other extended invariants which are less well known,
those of Krull dimension and Cohen-Macaulay defect.  The Krull
dimension can be found in both \cite{Foxby} and \cite{Iversen}, and
the Cohen-Macaulay defect in \cite{AvrFox}.  Let $M$ be in
$\Dfgp(\rk)$.  The Krull dimension of $M$ may be defined as
\begin{equation}
\label{equ:dim}
  \dim_{\rk} M = \sup \bigl( \dim_{\rk}\H_i(M) - i \bigr),
\end{equation}
see \cite[prop.\ 3.5]{Foxby}.  (Note the sign change induced by
the difference between the present homological notation and the
cohomological notation of \cite{Foxby}.)  The Cohen-Macaulay defect of
$M$ is then
\[
  \cmd_{\rk} M = \dim_{\rk} M - \depth_{\rk}M,
\]
see \cite[(1.1)]{AvrFox}.
\hfill $\Box$
\end{Notation}

\begin{Remark}
\label{rmk:generalization}
The final condition of Setup \ref{set:blanket} is that the canonical
ring homomorphism $\rk \rightarrow \H_0(R)$ is surjective.  

This implies that when localizing the Setup at a prime ideal of the
ground ring $\rk$, the ring $\H_0(R)$ remains local; a fact needed in
some of the proofs.

However, the results of the paper sometimes apply even if $\rk
\rightarrow \H_0(R)$ is not surjective. 

Namely, suppose that the conditions of Setup \ref{set:blanket} are
satisfied except that $\rk \rightarrow \H_0(R)$ is not surjective.
Suppose moreover that $R_0$ is central in $R$ and that $\H_0(R)$ is
finitely generated as an $\rk$-algebra by $\xi_1, \ldots, \xi_n$.  Then
the results of this paper still apply to $R$.

To see this, pick cycles $\Xi_1, \ldots, \Xi_n$ in $R_0$ representing
$\xi_1, \ldots, \xi_n$, set $\rk^{\prime} = \rk[X_1, \ldots, X_n]$,
and consider the $\rk$-linear ring homomorphism $\rk^{\prime}
\rightarrow R_0$ given by $X_i \mapsto \Xi_i$.  Then $R$ is a DGA over
$\rk^{\prime}$ and the canonical ring homomorphism $\rk^{\prime}
\rightarrow \H_0(R)$ is surjective.

To achieve the situation of Setup \ref{set:blanket}, it remains to
make the ground ring local.  For this, let $\fm$ be the maximal ideal
of $\H_0(R)$ and let $\fp$ be the contraction to $\rk^{\prime}$.
Replace $\rk^{\prime}$ and $R$ by the base changed versions
\[
  \widetilde{\rk} 
    = \rk^{\prime}_{\fp} \otimes_{\rk^{\prime}} \rk^{\prime} 
    \cong \rk^{\prime}_{\fp}
  \;\; \mbox{ and } \;\;
  \widetilde{R} 
    = \rk^{\prime}_{\fp} \otimes_{\rk^{\prime}} R.
\]
Then $\widetilde{\rk}$ is local and the canonical ring homomorphism
$\widetilde{\rk} \rightarrow \H_0(\widetilde{R})$ is surjective.
Hence the pair $\widetilde{\rk}$ and $\widetilde{R}$ fall under Setup
\ref{set:blanket}, and so the results of the paper apply to
$\widetilde{R}$. 

Now, the localization at $\fp$ inverts the elements of $\rk^{\prime}$
outside $\fp$.  Such elements are mapped to elements of $\H_0(R)$
outside $\fm$, and these are already invertible.  Hence the homology
of the canonical morphism
\[
  R \rightarrow \rk^{\prime}_{\fp} \otimes_{\rk^{\prime}} R 
\]
is an isomorphism; that is, the canonical morphism 
\[
  R \rightarrow \widetilde{R}
\]
is a qua\-si-i\-so\-mor\-phism.  This implies that $R$ and
$\widetilde{R}$ have equivalent derived categories, see
\cite[III.4.2]{KrizMay}, and so, since the results of this paper apply
to $\widetilde{R}$, they also apply to $R$.  \hfill $\Box$
\end{Remark}

\section{Homological estimates}
\label{sec:homological_estimates}

This Section provides some estimates which will be used as input for
the proofs of the main Theorems.

The following Lemma is well known.  It holds because $\H_0(R)$ is
local.  The proof is a simple application of the Eilenberg-Moore
spectral sequence, see \cite[exam.\ 1, p.\ 280]{FHTbook}.

\begin{Lemma}
\label{lem:inf}
Let $X$ be in $\Dfgp(R^{\opp})$ and let $Y$ be in $\Dfgp(R)$.  Then
\[
  \inf \bigl( X \LTensor_R Y \bigr) = \inf X + \inf Y.
\]
Consequently, if $X \not\cong 0$ and $Y \not\cong 0$ then $X
\LTensor_R Y \not\cong 0$.
\hfill $\Box$
\end{Lemma}

For the following Lemma, note that $\H_0(R)$ can be viewed as a DG
left-$\H_0(R)$-right-$R$-module; in subscript notation,
${}_{\H_0(R)}\H_0(R)_{R}$.  If $L = {}_{R}L$ is a DG left-$R$-module,
then
\[
  \H_0(R) \LTensor_R L = {}_{\H_0(R)}\H_0(R)_{R} \LTensor_R {}_{R}L
\]
inherits a DG left-$\H_0(R)$-module structure.  Since a DG
left-$\H_0(R)$-module is just a complex of left-$\H_0(R)$-modules,
$\H_0(R) \LTensor_R L$ is hence a complex of left-$\H_0(R)$-modules.

\begin{Lemma}
\label{lem:kpd}
Let $L$ be in $\D(R)$.  Then
\[
  L \mbox{ is in } \Dc(R) 
  \; \Leftrightarrow \;
  L \mbox{ is in } \Dfgp(R) \mbox{ and } \kpd_R L < \infty.
\]
If these equivalent statements hold, then $\H_0(R) \LTensor_R L$ is
in $\Dc(\H_0(R))$, and
\[
  \pd_{\H_0(R)} \bigl( \H_0(R) \LTensor_R L \bigr) = \kpd_R L.
\]
\end{Lemma}

\begin{proof}
$\Rightarrow\;$  Let $L$ be in $\Dc(R)$; that is, $L$ is finitely built
from $R$ in $\D(R)$.  Setup \ref{set:blanket} implies that $R$ is in
$\Dfgp(R)$.  Moreover, $\sup k \LTensor_R R = \sup k = 0 < \infty$.
But then $L$, being finitely built from $R$, is also in $\Dfgp(R)$ and
has $\sup k \LTensor_R L < \infty$.  And $\sup k \LTensor_R L <
\infty$ implies $\kpd_R L < \infty$ by \cite[rmk.\ (1.2)]{FJiic}. 

\medskip
\noindent
$\Leftarrow\;$  When $L$ is in $\Dfgp(R)$, there is a minimal semi-free
resolution $F \rightarrow L$ by \cite[(0.5)]{FJiic}.  When $\kpd_R L <
\infty$, it is not hard to see from \cite[(0.5) and lem.\
(1.7)]{FJiic} that there is a semi-free filtration of $F$ which only
contains finitely many quotients of the form $\Sigma^i R^{\alpha}$
where $\Sigma^i$ denotes the $i$'th suspension and where $\alpha$ is
finite.  This means that $F$ and hence $L$ is finitely built from $R$.

\medskip

Now suppose that the equivalent statements hold.  It is clear that
$\H_0(R) \LTensor_R R \cong \H_0(R)$ is in $\Dfgpm(\H_0(R))$.  As $L$
is finitely built from $R$, it follows that $\H_0(R) \LTensor_R L$ is
also in $\Dfgpm(\H_0(R))$.  Therefore the first $=$ in the following
computation holds by \cite[(A.5.7.3)]{LWC},
\begin{align*}
  \pd_{\H_0(R)} \bigl( \H_0(R) \LTensor_R L \bigr)
  & = - \inf \RHom_{\H_0(R)}(\H_0(R) \LTensor_R L,k) \\
  & = - \inf \RHom_R(L,k) \\
  & = \kpd_R L,
\end{align*}
where the second $=$ is by adjunction and the last $=$ is by
definition.
\end{proof}

The following Lemmas use that, as noted in Remark
\ref{rmk:A-structure}, all objects in sight have an $A$-structure.

\begin{Lemma}
\label{lem:depth}
Let $X$ be in $\Dfgm(R^{\opp})$ and let $L$ be in $\Dc(R)$.  Then
\[
  \depth_{\rk} \bigl( X \LTensor_R L \bigr) 
  =  \depth_{\rk} X - \kpd_R L.
\]
\end{Lemma}

\begin{proof}
The Lemma can be proved by a small variation of a well known proof of
the Auslander-Buchsbaum theorem, as given for instance in \cite[thm.\
3.2]{PJIdent}.  Let me give a summary for the benefit of the reader.

Since $L$ is finitely built from $R$ in $\D(R)$, there is an
isomorphism 
\[
  \RHom_{\rk}(k,X \LTensor_R L) 
  \cong \RHom_{\rk}(k,X) \LTensor_R L = (*).
\]
Replace $\RHom_{\rk}(k,X)$ with a quasi-isomorphic truncation $T$
concentrated in homological degrees $\leq \sup \RHom_{\rk}(k,X)$; see
\cite[(0.4)]{FJiic}.  Replace $L$ with a minimal semi-free resolution
$F$; see \cite[(0.5)]{FJiic}.  Then
\[
  (*) \cong T \otimes_R F,
\]
and hence
\[
  \sup \RHom_{\rk}(k,X \LTensor_R L) = \sup T \otimes_R F.
\]

The claim of the Lemma is that
\[
  \sup \RHom_{\rk}(k,X \LTensor_R L)
  = \sup \RHom_{\rk}(k,X) + \kpd_R L,
\]
and by the above this amounts to
\begin{equation}
\label{equ:j}
  \sup T \otimes_R F = \sup \RHom_R(k,X) + \kpd_R L.
\end{equation}

Forgetting the differentials of $R$ and $F$ gives the underlying
graded algebra $R^{\natural}$ and the underlying graded module
$F^{\natural}$, and \cite[(0.5)]{FJiic} says that
\[
  F^{\natural} \cong 
  \coprod_{i \leq \operatorname{k.pd}_R L} \Sigma^i(R^{\natural})^{\beta_i}.
\]
Hence
\[
  (T \otimes_R F)^{\natural}
  \cong T^{\natural} \otimes_{R^{\natural}} F^{\natural}
  \cong T^{\natural} \otimes_{R^{\natural}} \coprod_{i \leq \operatorname{k.pd}_R L} \Sigma^i(R^{\natural})^{\beta_i}
  \cong \coprod_{i \leq \operatorname{k.pd}_R L} \Sigma^i(T^{\natural})^{\beta_i}.
\]

The right hand side is just a collection of copies of $T^{\natural}$
moved around by $\Sigma^i$, so since $T$ and hence $T^{\natural}$ is
concentrated in homological degrees $\leq \sup \RHom_{\rk}(k,X)$, the
right hand side and therefore the left hand side is concentrated in
homological degrees $\leq \sup \RHom_{\rk}(k,X) + \kpd_R L$.  This
implies 
\[
  \sup T \otimes_R F \leq \sup \RHom_{\rk}(k,X) + \kpd_R L.
\]

Using that $\beta_{\operatorname{k.pd}_R L} \not= 0$ by \cite[lem.\
(1.7)]{FJiic}, it is possible also to see
\[
  \sup T \otimes_R F \geq \sup \RHom_{\rk}(k,X) + \kpd_R L.
\]
This proves Equation \eqref{equ:j} and hence the Lemma.
\end{proof}

Through the canonical morphism $R \rightarrow \H_0(R)$, an
$\H_0(R)$-module $M$ can be viewed as a DG right-$R$-module.  If $M$
is finitely generated over $\H_0(R)$, then as a DG right-$R$-module it
is in $\Dfgpm(R^{\opp})$. 

\begin{Lemma}
\label{lem:dim}
Let $M$ be a finitely generated $\H_0(R)$-module and let $L$ be in
$\Dfgp(R)$.  Suppose $M \not\cong 0$ and $L \not\cong 0$.  View $M$ as a DG
right-$R$-module in $\Dfgpm(R^{\opp})$, and suppose
\[
  \dim_{\rk} \H_i \bigl( M \LTensor_R L \bigr) \leq 0
\]
for each $i$.  Then
\[
  \kpd_R L  \geq  \dim_{\rk}M + \inf L.
\]
\end{Lemma}

\begin{proof} 
If $\kpd_R L = \infty$ then the Lemma holds trivially, so suppose
$\kpd_R L < \infty$.  Then Lemma \ref{lem:kpd} says that $\H_0(R)
\LTensor_R L$ is in $\Dc(\H_0(R))$.  That is, $\H_0(R) \LTensor_R L$ is
finitely built from $\H_0(R)$, so $\H_0(R) \LTensor_R L$ is isomorphic
to a bounded complex of finitely generated free $\H_0(R)$-modules.
Also, Lemma \ref{lem:inf} implies $\H_0(R) \LTensor_R L \not\cong 0$, and
hence \cite[thm.\ 4.1]{Iversen} says
\begin{align*}
  \pd_{\H_0(R)} \bigl( \H_0(R) \LTensor_R L \bigr) \geq
  & \dim_{\H_0(R)} M \\
\label{equ:a}
  & - \dim_{\H_0(R)} \bigl( M \LTensor_{\H_0(R)} (\H_0(R) \LTensor_R L) \bigr).
\end{align*}
Note that the assumption in \cite{Iversen} that the ring is
equicharacteristic is unnecessary: The assumption is only used to
ensure that the so-called new intersection theorem is valid, and this
was later proved for all local noetherian commutative rings in
\cite[thm.\ 1]{Roberts}.

Moving the parentheses in the last term gets rid of tensoring with
$\H_0(R)$, and Krull dimensions over $\H_0(R)$ can be replaced with
Krull dimensions over $\rk$ because $\rk \rightarrow \H_0(R)$ is
surjective, so the inequality is
\begin{equation}
\label{equ:a}
  \pd_{\H_0(R)} \bigl( \H_0(R) \LTensor_R L \bigr) \geq
  \dim_{\rk} M
  - \dim_{\rk} \bigl( M \LTensor_R L \bigr).
\end{equation}
The first term here is
\begin{equation}
\label{equ:b}
  \pd_{\H_0(R)} \bigl( \H_0(R) \LTensor_R L \bigr) = \kpd_R L
\end{equation}
by Lemma \ref{lem:kpd}.  For the third term, note that the assumption
\[
  \dim_{\rk} \H_i \bigl( M \LTensor_R L \bigr) \leq 0
\]
for each $i$ implies
\begin{equation}
\label{equ:c}
  \dim_{\rk} \bigl( M \LTensor_R L \bigr) 
  = - \inf \bigl( M \LTensor_R L \bigr) = (*);
\end{equation}
see Notation \ref{not:dim}.  But Lemma \ref{lem:inf} implies
\begin{equation}
\label{equ:d}
  (*) = - \inf L.
\end{equation}

Substituting Equations \eqref{equ:b} to \eqref{equ:d} into the
inequality \eqref{equ:a} gives the inequality claimed in the Lemma.
\end{proof}

\section{Flat base change}
\label{sec:flat_base_change}

This Section sets up a theory of flat base change which will be used
in the proofs of the main Theorems.

Let $\widetilde{\rk}$ be a local noetherian commutative ring and let
$\rk \rightarrow \widetilde{\rk}$ be a flat ring homomorphism.

It is clear that
\[
  \widetilde{R} = \widetilde{\rk} \otimes_{\rk} R
\]
is a chain DGA over $\widetilde{\rk}$.  The homology is
\[
  \H_i(\widetilde{R})
  = \H_i \bigl( \widetilde{\rk} \otimes_{\rk} R \bigr) 
  \cong \widetilde{\rk} \otimes_{\rk} \H_i(R)
\]
and this is finitely generated over $\widetilde{\rk}$ for each $i$.
The canonical ring homomorphism $\rk \rightarrow \H_0(R)$ is
surjective, so $\widetilde{\rk} \otimes_{\rk} \rk \rightarrow
\widetilde{\rk} \otimes_{\rk} \H_0(R)$ is also surjective, but this
map is isomorphic to the canonical ring homomorphism
\[
  \widetilde{\rk} \rightarrow \H_0(\widetilde{R})
\]
which is hence surjective.  So Setup \ref{set:blanket} applies to the
DGA $\widetilde{R}$ over the ring $\widetilde{\rk}$.

There is a morphism of DGAs 
\[
  R \rightarrow \widetilde{R}
\]
given by $r \mapsto 1 \otimes r$, and this defines a base change
functor of DG left modules
\[
  \widetilde{R} \LTensor_R - 
  : \D(R) \rightarrow \D(\widetilde{R})
\]
which in fact is just given by
\begin{equation}
\label{equ:base_change_functors}
  \widetilde{R} \LTensor_R -
  = (\widetilde{\rk} \otimes_{\rk} R) \LTensor_R -
  \cong \widetilde{\rk} \otimes_{\rk} (R \LTensor_R -)
  \cong \widetilde{\rk} \otimes_{\rk} -.
\end{equation}
There is also a base change functor of DG right modules.  

It is easy to see that the base change functors preserve membership of
the subcategories $\Dc$, $\Dfgp$, $\Dfgm$, and $\Dfgpm$.

If $X$ is in $\D(R^{\opp})$ and $Y$ is in $\D(R)$, then it is an
exercise to compute the derived tensor product of the base changed DG
modules $\widetilde{X} = X \LTensor_R \widetilde{R}$ and
$\widetilde{Y} = \widetilde{R} \LTensor_R Y$ as
\begin{equation} 
\label{equ:tensor}
  \widetilde{X} \LTensor_{\widetilde{R}} \widetilde{Y}
  \cong \widetilde{\rk} \otimes_{\rk} (X \LTensor_R Y).
\end{equation}

\section{Amplitude inequalities for compact objects}
\label{sec:amp_ineq_A}

This Section proves Theorem \ref{thm:amp} which is a DGA
generalization of the first of Iversen's amplitude inequalities from
\cite{Iversen}.  Theorems A and A${}^{\prime}$ from the Introduction
follow easily.

\begin{Theorem}
\label{thm:amp}
Let $X$ be in $\Dfgpm(R^{\opp})$ and let $L$ be in $\Dc(R)$.  Suppose $X
\not\cong 0$ and $L \not\cong 0$.  Then
\[
  \amp \bigl( X \LTensor_R L \bigr) \geq \amp X.
\]
\end{Theorem}

\begin{proof}
The inequality says
\[
  \sup \bigl( X \LTensor_R L \bigr) - \inf \bigl( X \LTensor_R L \bigr)
  \geq
  \sup X - \inf X,
\]
which by Lemma \ref{lem:inf} is the same as
\begin{equation}
\label{equ:e}
  \sup \bigl( X \LTensor_R L \bigr) \geq \sup X + \inf L.
\end{equation}

Write
\[
  M = \H_{\sup X}(X)
\]
for the top homology of $X$.  With this notation, \cite[prop.\
3.17]{Foxby} says
\begin{equation}
\label{equ:f}
  \depth_{\rk} X  \leq  \dim_{\rk} M - \sup X.
\end{equation}
(Note again the difference between the homological notation of this
paper and the cohomological notation of \cite{Foxby}.)

To prove the Theorem, consider first the special case where
\[
  \dim_{\rk} \H_i \bigl( M \LTensor_R L \bigr) \leq 0
\]
for each $i$.  Then
\begin{align*}
  \sup \bigl( X \LTensor_R L \bigr)
  & \stackrel{\rm (a)}{\geq}
      - \depth_{\rk} \bigl( X \LTensor_R L \bigr) \\
  & \stackrel{\rm (b)}{=}
      - \depth_{\rk} X + \kpd_R L \\
  & \stackrel{\rm (c)}{\geq}
      - \depth_{\rk} X + \dim_{\rk} M + \inf L \\
  & \stackrel{\rm (d)}{\geq}
      - \dim_{\rk} M + \sup X + \dim_{\rk} M + \inf L \\
  & = \sup X + \inf L
\end{align*}
proving \eqref{equ:e}.  Here (a) is by \cite[eq.\ (3.3)]{Foxby}, (b)
is by Lemma \ref{lem:depth}, (c) is by Lemma \ref{lem:dim}, and (d) is
by Equation \eqref{equ:f}.

Next the general case which will be reduced to the above special case
by localization.  Observe that $M \LTensor_R L \not\cong 0$ by Lemma
\ref{lem:inf}.  Pick a prime ideal $\fp$ of $\rk$ which is minimal in
\begin{equation}
\label{equ:g}
  \bigcup_i \Supp_{\rk} \H_i \bigl( M \LTensor_R L \bigr)
\end{equation}
and consider the flat ring homomorphism $\rk \rightarrow \rk_{\fp}$.
Set 
\[
  \widetilde{R} = \rk_{\fp} \otimes_{\rk} R, \;\;\;\;\;
  \widetilde{X} = \rk_{\fp} \otimes_{\rk} X, \;\;\;\;\;
  \widetilde{L} = \rk_{\fp} \otimes_{\rk} L
\]
so that $\widetilde{X}$ and $\widetilde{L}$ are the base changes of $X$
and $L$ to $\widetilde{R}$, see Section \ref{sec:flat_base_change}.

Let me check that the above special case of the Theorem applies to
$\widetilde{X}$ and $\widetilde{L}$.  The theory of Section
\ref{sec:flat_base_change} says that Setup \ref{set:blanket} applies
to $\widetilde{R}$ over $\rk_{\fp}$, that $\widetilde{X}$ is in
$\Dfgpm(\widetilde{R}^{\opp})$, and that $\widetilde{L}$ is in
$\Dc(\widetilde{R})$.  Moreover, $\fp$ is in the support of some $\H_i
\bigl( M \LTensor_R L \bigr)$ in $\rk$ so must be in the support of $M
= \H_{\sup X}(X)$ and in the support of some $\H_i(L)$.  It follows
that $\widetilde{X} \not\cong 0$ and $\widetilde{L}
\not\cong 0$.

Since $\fp$ is in the support of $M = \H_{\sup X}(X)$, it even follows
that 
\begin{equation}
\label{equ:h}
  \sup \widetilde{X} = \sup X
\end{equation}
and
\[
  \widetilde{M} = \H_{\sup \widetilde{X}}(\widetilde{X})
  = \H_{\sup X} \bigl( \rk_{\fp} \otimes_{\rk} X \bigr)
  \cong \rk_{\fp} \otimes_{\rk} \H_{\sup X}(X)
  = \rk_{\fp} \otimes_{\rk} M.
\]
Finally, Equation \eqref{equ:tensor} from Section
\ref{sec:flat_base_change} implies
\[
  \H_i \bigl( \widetilde{M} \LTensor_{\widetilde{R}} \widetilde{L} \bigr)
  \cong \rk_{\fp} \otimes_{\rk} \H_i \bigl( M \LTensor_R L \bigr).
\]
The support of each of these modules in $\rk_{\fp}$ is either empty or
equal to the maximal ideal $\fp_{\fp}$ since $\fp$ was chosen minimal
in the set \eqref{equ:g}, and each of the modules is finitely
generated over $\rk_{\fp}$ because each $\H_i(M \LTensor_R L)$ is
finitely generated over $\H_0(R)$ and hence over $\rk$.  So
\[
  \dim_{\rk_{\fp}}
  \H_i \bigl( \widetilde{M} \LTensor_{\widetilde{R}} \widetilde{L} \bigr) \leq 0
\]
for each $i$.  

Hence the above special case of the Theorem does apply and gives
\[
  \sup \bigl( \widetilde{X} \LTensor_{\widetilde{R}} \widetilde{L} \bigr)
  \geq \sup \widetilde{X} + \inf \widetilde{L},
\]
which by Equation \eqref{equ:tensor} again is
\begin{equation}
\label{equ:i}
  \sup \bigl( \rk_{\fp} \otimes_{\rk} (X \LTensor_R L) \bigr)
  \geq \sup \widetilde{X} + \inf \widetilde{L}.
\end{equation}
So
\begin{align*}
  \sup \bigl( X \LTensor_R L \bigr)
  & \geq \sup \bigl( \rk_{\fp} \otimes_{\rk} (X \LTensor_R L) \bigr) \\
  & \stackrel{\rm (e)}{\geq}
      \sup \widetilde{X} + \inf \widetilde{L} \\
  & \stackrel{\rm (f)}{=}
      \sup X + \inf \widetilde{L} \\
  & = \sup X + \inf \bigl( \rk_{\fp} \otimes_{\rk} L \bigr) \\
  & \geq \sup X + \inf L
\end{align*}
proving \eqref{equ:e}.  Here (e) is by \eqref{equ:i} and (f) is by
\eqref{equ:h}. 
\end{proof}

{
\noindent
{\it Proof } (of Theorems A and A${}^{\prime}$).
Theorem A follows by setting $X = R$ in Theorem \ref{thm:amp}, and
Theorem A${}^{\prime}$ is equivalent to Theorem A by Lemma
\ref{lem:kpd}. 
\hfill $\Box$
\medskip
}

\section{Amplitude inequality for objects with finite $k$-injective
dimension}
\label{sec:amp_ineq_B}

This Section proves Theorem B from the Introduction.  The proof uses
dualizing complexes; see \cite[chp.\ V]{Hartshorne}.  Since, on one
hand, not all rings have dualizing complexes, while, on the other,
complete local noetherian commutative rings do, it is also necessary
to include some material on completions.

The following Proposition assumes that the ground ring $\rk$ has a
dualizing complex $C$, and considers the DG left-$R$-right-$R$-module 
\[
  {}_{R}D_{R} = \RHom_{\rk}({}_{R}R_{R},C)
\]
whose left-structure comes from the right-structure of the $R$ in the
first argument of $\RHom$, and vice versa.  By forgetting the
right-structure, I can get a DG left-$R$-module ${}_{R}D$.

\begin{Proposition}
\label{pro:kid}
Suppose that $\rk$ has a dualizing complex $C$ and set ${}_{R}D_{R} =
\RHom_{\rk}({}_{R}R_{R},C)$.  Let $I$ be in $\Dfgm(R)$.  Then the
following conditions are equivalent.
\begin{enumerate}

\smallskip
  \item  $\kid_R I < \infty$.

\smallskip
  \item  ${}_{R}I$ is finitely built from ${}_{R}D$ in $\D(R)$.

\smallskip
  \item  ${}_{R}I \cong {}_{R}D_{R} \LTensor_R {}_{R}L$
         for an ${}_{R}L$ in $\Dc(R)$.

\end{enumerate}
\end{Proposition}

\begin{proof}
(i) $\Rightarrow$ (iii).  Let
\[
  M_R = \RHom_{\rk}({}_{R}I,C)
\]
be the dual of $I$.  Since $I$ is in $\Dfgm(R)$ and hence in
$\Dfgm(\rk)$, it follows that $M$ is $\Dfgp(\rk)$ and hence in
$\Dfgp(R^{\opp})$.  Moreover,
\[
  M \LTensor_R k
  = \RHom_{\rk}(I,C) \LTensor_R k
  \stackrel{\rm (a)}{\cong} \RHom_{\rk}(\RHom_R(k,I),C) = (*),
\]
where (a) holds because $C$ is isomorphic in $\D(\rk)$ to a bounded
complex of injective modules, cf.\ \cite[(A.4.24)]{LWC}.  The
assumption $\kid_R I < \infty$ implies that $\RHom_R(k,I)$ has bounded
homology, so the same is true for $(*)$ whence $\sup M
\LTensor_R k < \infty$.  This implies $\kpd_{R^{\opp}} M < \infty$ by
\cite[rmk.\ (1.2)]{FJiic}.

By Lemma \ref{lem:kpd} this means that $M$ is in $\Dc(R^{\opp})$, that
is, $M_R$ is finitely built from $R_R$.  But then ${}_{R}L =
\RHom_{R^{\opp}}(M_{R},{}_{R}R_{R})$ is finitely built from ${}_{R}R$
and satisfies
\[
  M_{R} \cong \RHom_R({}_{R}L,{}_{R}R_{R}), 
\]
and hence 
\begin{align*}
  {}_{R}I & \stackrel{\rm(b)}{\cong} \RHom_{\rk}(\RHom_{\rk}({}_{R}I,C),C) \\
          & = \RHom_{\rk}(M_R,C) \\
          & \cong \RHom_{\rk}(\RHom_R({}_{R}L,{}_{R}R_{R}),C) \\
          & \stackrel{\rm (c)}{\cong} \RHom_{\rk}({}_{R}R_{R},C) \LTensor_R {}_{R}L \\
          & = {}_{R}D_{R} \LTensor_R {}_{R}L,
\end{align*}
proving (iii).  Here (b) is by \cite[thm.\ (A.8.5)]{LWC} and (c) is
because ${}_{R}L$ is finitely built from ${}_{R}R$.

\medskip
\noindent
(iii) $\Rightarrow$ (ii).  For ${}_{R}L$ to be in $\Dc(R)$ means that
${}_{R}L$ is finitely built from ${}_{R}R$ in $\D(R)$.  But then
${}_{R}I \cong {}_{R}D_{R} \LTensor_R {}_{R}L$ is finitely built from 
\[
  {}_{R}D_{R} \LTensor_R {}_{R}R \cong {}_{R}D
\]
in $\D(R)$.

\medskip
\noindent
(ii) $\Rightarrow$ (i).  Without loss of generality, I can assume that
$C$ is nor\-ma\-li\-zed, that is, $\RHom_{\rk}(k,C) \cong k$.  Then
\begin{align*}
  \RHom_R(k,{}_{R}D) 
  & = \RHom_R(k,\RHom_{\rk}(R_{R},C)) \\
  & \stackrel{\rm (d)}{\cong} \RHom_{\rk}(R_{R} \LTensor_R k,C) \\
  & \cong \RHom_{\rk}(k,C) \\
  & \cong k
\end{align*}
has bounded homology, where (d) is by adjunction.  When ${}_{R}I$ is
finitely built from ${}_{R}D$, then the homology of
$\RHom_R(k,{}_{R}I)$ is also bounded, and then
\[
  \kid_R I = - \inf \RHom_R(k,{}_{R}I) < \infty.
\]
\end{proof}

For the remaining part of the paper, let $\fm$ be the maximal ideal of
$\rk$ and consider $\widehat{\rk}$, the completion of $\rk$ in the
$\fm$-adic topology, which is a local noetherian commutative ring by
\cite[p.\ 63, (4)]{Matsumura}.

The canonical ring homomorphism $\rk \rightarrow \widehat{\rk}$ is
flat by \cite[p.\ 63, (3)]{Matsumura}, and the theory of Section
\ref{sec:flat_base_change} gives a new chain DGA
\[
  \widehat{R} = \widehat{\rk} \otimes_{\rk} R
\]
over $\widehat{\rk}$, and base change functors for DG modules from $R$
to $\widehat{R}$.

\begin{Lemma}
\label{lem:completion}
Let $I$ be in $\D(R)$ and consider the base changed DG module
$\widehat{I} = \widehat{R} \LTensor_R I$ in $\D(\widehat{R})$.  Then 
\begin{enumerate}

  \item  $\amp \widehat{I} = \amp I$.

\smallskip
  \item  $\kid_{\widehat{R}} \widehat{I} = \kid_R I$.

\end{enumerate}
\end{Lemma}

\begin{proof}
(i).  Equation \eqref{equ:base_change_functors} from Section
\ref{sec:flat_base_change} says
\[
  \widehat{I} = \widehat{R} \LTensor_R I 
  \cong \widehat{\rk} \otimes_{\rk} I.
\]
Since $\widehat{\rk}$ is faithfully flat over $\rk$ by \cite[p.\ 63,
(3)]{Matsumura}, part (i) is clear.

\medskip
\noindent
(ii).  The residue class field of $\H_0(R)$ is $k$, so the residue
class field of $\H_0(\widehat{R}) \cong \widehat{\rk} \otimes_{\rk}
\H_0(R)$ is $\widehat{\rk} \LTensor_{\rk} k$ which by Equation
\eqref{equ:base_change_functors} is
\[
  \widehat{\rk} \otimes_{\rk} k \cong \widehat{R} \LTensor_R k.
\]
Hence
\[
  \kid_{\widehat{R}} \widehat{I}
  = -\inf \RHom_{\widehat{R}}(\widehat{R} \LTensor_{R} k,\widehat{I})
  = (*).
\]

But
\begin{align*}
  \RHom_{\widehat{R}}(\widehat{R} \LTensor_R k,\widehat{I})
  & \stackrel{\rm (a)}{\cong} \RHom_R(k,\RHom_{\widehat{R}}(\widehat{R},\widehat{I})) \\
  & \cong \RHom_R(k,\widehat{I}) \\
  & \cong \RHom_R(k,\widehat{\rk} \otimes_{\rk} I) \\
  & \stackrel{\rm (b)}{\cong}
    \widehat{\rk} \otimes_{\rk} \RHom_R(k,I),
\end{align*}
where (a) is by adjunction and (b) is because $\widehat{\rk}$ is flat
over $\rk$ while $k$ has finitely generated homology, cf.\
\cite[(A.4.23)]{LWC}.  Hence
\[
  (*) = - \inf \widehat{\rk} \otimes_{\rk} \RHom_R(k,I)
  \stackrel{\rm (c)}{=} - \inf \RHom_R(k,I)
  = \kid_R I,
\]
proving part (ii).  Here (c) is again because $\widehat{\rk}$ is
faithfully flat over $\rk$.
\end{proof}

{
\noindent
{\it Proof } (of Theorem B).
The base change $\rk \rightarrow \widehat{\rk}$ induces the change
from $R$ and $I$ to $\widehat{R}$ and $\widehat{I}$.  

$\widehat{I}$ is in $\Dfgm(\widehat{R})$ by the theory of Section
\ref{sec:flat_base_change}.  Lemma \ref{lem:completion} implies that
$\kid_{\widehat{R}}\widehat{I} < \infty$ and that $\amp
\widehat{I} = \amp I$.  Moreover, $\cmd_{\widehat{\rk}}\widehat{R} =
\cmd_{\rk}R$ by \cite[prop.\ (1.2)]{AvrFox}.  So it is enough to prove
Theorem B for $\widehat{I}$ over $\widehat{R}$.

Setup \ref{set:blanket} applies to $\widehat{R}$ over $\widehat{\rk}$
by Section \ref{sec:flat_base_change}, so the results proved so far
apply to DG modules over $\widehat{R}$.  Since $\widehat{\rk}$ is
complete, it has a dualizing complex $C$ by \cite[sec.\
V.10.4]{Hartshorne}.  Hence Proposition \ref{pro:kid} gives
\[
  \widehat{I} \cong 
  \RHom_{\widehat{\rk}}(\widehat{R},C) \LTensor_{\widehat{R}} L
\]
for an $L$ in $\Dc(\widehat{R})$.  But then Theorem \ref{thm:amp}
gives 
\[
  \amp \widehat{I} 
  \geq \amp \RHom_{\widehat{\rk}}(\widehat{R},C)
  = \cmd_{\widehat{\rk}} \widehat{R}
\]
as desired, where the $=$ is by \cite[(1.3.2)]{AvrFox}.
\hfill $\Box$
\medskip
}

\end{document}